\documentclass[12pt,a4paper]{amsart}
\usepackage{amssymb,amsmath,amsthm}
\usepackage{graphicx}
\usepackage{hyperref}
 \usepackage[color=green!40]{todonotes}

\newcommand\cF{{\mathcal F}}
\newcommand\cS{{\mathcal S}}
\newcommand\x{{\mathbf x}}
\newcommand\y{{\mathbf y}}
\newcommand\z{{\mathbf z}}

\newcommand\cG{{\mathcal G}}

\makeatletter
\newtheorem*{rep@theorem}{\rep@title}
\newcommand{\newreptheorem}[2]{%
\newenvironment{rep#1}[1]{%
 \def\rep@title{#2 \ref{##1}}%
 \begin{rep@theorem}}%
 {\end{rep@theorem}}}
\makeatother

\theoremstyle{plain}
\newtheorem{theorem}{Theorem}[section]
\newreptheorem{theorem}{Theorem}
\newtheorem{lemma}[theorem]{Lemma}
\newtheorem{corollary}[theorem]{Corollary}

\newtheorem{proposition}[theorem]{Proposition}
\newtheorem{obs}[theorem]{Observation}

\theoremstyle{definition}

\newtheorem{defn}[theorem]{Definition}
\newtheorem{con}[theorem]{Construction}

\newcommand\cref[1]{Corollary~\ref{cor:#1}}

\title{Vector sum-intersection theorems}

\author{Bal\'azs Patk\'os}
\address{Alfr\'ed R\'enyi Institute of Mathematics}
\email{patkos@renyi.hu}
\thanks{Patk\'os's research is partially supported by NKFIH grants SNN 129364 and FK 132060.}

\author{Zsolt Tuza}
\address{Alfr\'ed R\'enyi Institute of Mathematics and University of Pann\'onia}
\email{tuza.zsolt@mik.uni-pannon.hu}
\thanks{Tuza's research is partially supported by NKFIH grant SNN 129364.}

\author{M\'at\'e Vizer}
\address{Alfr\'ed R\'enyi Institute of Mathematics}
\email{vizermate@gmail.com}
\thanks{Vizer's research is partially supported by NKFIH grants SNN 129364, FK 132060, KH130371, by the  J\'anos Bolyai Research Fellowship  and by the New National Excellence Program under the grant number \'UNKP-21-5-BME-361.}

\date{}
\begin{document}

\maketitle

\begin{abstract}
     We introduce the following generalization of set intersection via characteristic vectors: for $n,q,s, t \ge 1$ a family $\cF\subseteq \{0,1,\dots,q\}^n$ of vectors is said to be \emph{$s$-sum $t$-intersecting} if for any distinct $\x,\y\in \cF$ there exist at least $t$ coordinates, where the entries of $\x$ and $\y$ sum up to at least $s$, i.e.\ $|\{i:x_i+y_i\ge s\}|\ge t$.
     The original set intersection corresponds to the case $q=1,s=2$. 
     
     We address analogs of several variants of classical results in this setting: the Erd\H os--Ko--Rado theorem and the theorem of Bollob\'as on intersecting set pairs.
     
\end{abstract}

\section{Introduction}

Many problems in extremal finite set theory ask for the maximum size of a set family that satisfies some intersection property. When members of the family examined are subsets of $[n]:=\{1,2,\dots,n\}$, then there is a one-to-one correspondence between a set $F$ and its $0$-$1$ characteristic vector $\x_F$ of length $n$, that has a 1-entry in its $i$th coordinate if and only if $i\in F$ for $i \in [n]$. 
So one can say that two sets $F$ and $G$ intersect, if the sum of their characteristic vectors (as vectors in $\mathbb{Z}^n$) contains a $2$ in at least one coordinate. 
The goal of this paper is to introduce a notion of intersection that generalizes set intersection (translated to sum of characteristic vectors) to a type of intersection among $q$-ary vectors.

For $q,n \ge 1$ we introduce the notation $Q=Q(q)=\{0,1,\dots,q\}$ and also $Q^n:=\{0,1,\dots,q\}^n$. 
We denote vectors by boldface letters and the $i$th coordinate of the vector $\x$ is denoted by $x_i$.

Intersection problems have been studied for vectors / integer sequences with several possible definitions for the size of the intersection: the \textit{permutation-type} intersection size of $\x,\y\in \{0,1,\dots,q\}^n$ is $|\x \cap_{\textrm{perm}} \y|=|\{i:x_i=y_i\}|$; the \textit{multiset-type} intersection size is defined as $|\x \cap_{\textrm{multi}} \y|=\sum_i \min\{x_i,y_i\}$. Results using the former definition include \cite{ff80,ft99}, while multiset-type results can be found in e.g.\ \cite{ft16,fgv}. The main definition of our paper is as follows.



\begin{defn}\label{ssumintersect}
For integers $n,q, s \ge 1$ and two vectors $\x,\y\in Q^n$, we define the \textit{size of their $s$-sum intersection} as $|\x\cap_s \y|=|\{i:x_i+y_i\ge s\}|$.


For $t \ge 1$ we say that $\x, \y \in Q^n$ are \textit{$s$-sum $t$-intersecting}, if $|\x\cap_s \y| \ge t$. More generally, $\cF \subset Q^n$ is \textit{s-sum t-intersecting} if any two vectors $\x, \y \in \cF$ are $s$-sum $t$-intersecting. 

In case of $t=1$ we just simply write \textit{$s$-sum intersecting} instead of $s$-sum $1$-intersecting.

\end{defn}

Note that the case $q=1$ and $s=2$ corrseponds to ordinary set intersection. 


\ 

We will consider analogs of the Erdős--Ko--Rado theorem and theorems about Bollobás's intersecting set-pair systems. To be able to state our results first we need to define uniformity for families of vectors. One has several options: as in the case of multisets and many other types of problems, we can work with the \textit{weight/rank} $r(\x):=\sum_{i=1}^n x_i$ of $\x \in Q^n$ and say that for an integer $r \ge 0$ a family $\cF\subseteq Q^n$ is \textit{$r$-rank uniform} if $r(\x)=r$ for all $\x\in \cF$. Another possibility is to use the size of the \textit{support} $S_\x=\{i:x_i\neq 0\}$ of $\x$. We say that $\cF \subseteq Q^n$ is \textit{$r$-support uniform} if $|S_\x|=r$ for every $\x \in \cF$.

\medskip

We use standard notation. 
For any set $X$, we denote by $\binom{X}{r}$ the family of all $r$-subsets of $X$ and $2^X$ denotes the power set of $X$.
For a set $F \subset [n]$ we denote its complement $[n] \setminus F$ by $\overline{F}$ and for $\cF$ a family of subsets of $[n]$ we introduce the notation $\overline{\cF}:=\{\overline{F} : F \in \cF \}$. 

As a vector analog for any  $\x \in Q^n$ we define its `complement' $\overline{\x}$ by letting $\overline{x}_i:=q-x_i$ for all $i \in [n]$ and for  a family $\cF$ of vectors in $Q^n$ we write $\overline{\cF}:=\{ \overline{\x} : \x \in \cF \}$. 






\medskip

The structure of the paper is as follows. In Subsection 1.1 we state various results about $s$-sum intersecting families of vectors, while in Subsection 1.2 we list our results about intersecting vector pairs. 
In Section 2 and Section 3 we prove our results about intersecting vectors and intersecting vector pairs, respectively. 
In Section 4---as concluding results---we give a new intersection definition to provide analogs of some results that would not work with $s$-sum intersection.

\subsection{Results on intersecting families of vectors}

Let us start with stating the seminal result of Erdős, Ko and Rado \cite{ekr61}.  


\begin{theorem}[Erd\H os, Ko, Rado \cite{ekr61}]\label{ekr}
For $n, r \ge 1$ with $2r\le n$ if $\cF\subseteq \binom{[n]}{r}$ is an intersecting family, then $|\cF|\le \binom{n-1}{r-1}$. Moreover, if $2r<n$ and $|\cF|=\binom{n-1}{r-1}$, then $\cF=\cF_i:=\{F:i\in F\in \binom{[n]}{r}\}$ holds for some $i\in [n]$. 

Furthermore, for any $1\le t<r$ there exists $n_0=n_0(r,t)$ such that if $\cF\subseteq \binom{[n]}{r}$ is $t$-intersecting, then $|\cF|\le \binom{n-t}{r-t}$, and equality holds if and only if $\cF=\{F: T\subset F\in \binom{[n]}{r}\}$ for some $T\in \binom{[n]}{t}$.
\end{theorem}

The exact value of the smallest possible $n_0(r,t)$ was obtained by Frankl \cite{f78} and Wilson~\cite{w84}. The largest possible size of an $r$-uniform $t$-intersecting family for all values of $n,t,r \ge 1$ was determined by Ahlswede and Khachatrian \cite{ak}.

Our first result is a generalization of the Erdős--Ko--Rado (EKR) theorem for $r$-support uniform families. If $s$ is even, then the vector family corresponding to $\cF_i$ of Theorem \ref{ekr} is $\{\x\in Q^n:x_i\ge \frac{s}{2}, |S_\x|=r\}$. If $s$ is odd, then to the family $\{\x\in Q^n:x_i\ge \lceil\frac{s}{2}\rceil, |S_\x|=r\}$ one can add vectors $\y$ with $y_i=\lfloor \frac{s}{2}\rfloor$ that pairwise $s$-sum intersect on some other coordinate. 

Observe that if $s\le q$ holds, then any vector having at least one entry at least $s$ can be added to any $s$-sum intersecting family, so it is enough to consider the case $q<s$.

\begin{theorem}\label{qekr}
For any $2\le q<s$ and integer $r \ge 1$,
if $\cF\subseteq Q^n$ is $r$-support uniform $s$-sum intersecting with $n\ge qr^2$, then 
\begin{eqnarray}
|\cF|\le\left\{
\begin{array}{cc} 
(q-\frac{s}{2}+1) q^{r-1}\binom{n-1}{r-1} & \textnormal{if} ~s ~\textnormal{is even,} \\
1+(q-\lceil\frac{s}{2}\rceil+1)\sum_{i=1}^r \binom{n-i}{r-i}q^{r-i}& \textnormal{if $s$ is odd,}
\end{array}
\right.
\end{eqnarray}
and these bounds are best possible.
\end{theorem}

The statement and proof of Theorem \ref{qekr} can be adjusted for the $r$-rank uniform case, too. 
Instead, we provide a different proof in the special case $s=q+1$ that works for all meaningful values of $n$.
Before stating our theorem, observe that if both $\x,\y\in Q^n$ have rank less than $\frac{q+1}{2}$, then they cannot $(q+1)$-sum intersect, while if both of them have rank greater than $\frac{qn}{2}$, then they always $(q+1)$-sum intersect. We denote by $Q(n,r)$ the set of all vectors in $Q^n$ of rank $r$. Just as in Theorem \ref{qekr}, if $q+1$ is even, then extremal families are stars (all vectors with entry at least $\frac{q+1}{2}$ in one fixed coordinate), while if $q+1$ is odd, then one can add further vectors to the star.

\begin{theorem}\label{rankqekr}
Let $n,q,r \ge 1$ and $\cF\subseteq Q^n$ be an $r$-rank uniform $(q+1)$-sum intersecting family with $\frac{q+1}{2}\le r\le \frac{qn}{2}$. Then
\begin{eqnarray}
|\cF|\le\left\{
\begin{array}{cc} 
\sum_{j=\frac{q+1}{2}}^q|Q(n-1,r-j)| & \textnormal{if} ~q+1 ~\textnormal{is even}, \\
1+\sum_{j=\lceil\frac{q+1}{2}\rceil}^q\sum_{i=1}^{\lfloor \frac{2(r-1)}{q}\rfloor}|Q(n-i,r-j-\frac{(i-1)q}{2})|& \textnormal{if $q+1$ is odd},
\end{array}
\right.
\end{eqnarray}
and these bounds are best possible.
\end{theorem}

Now we continue with $s$-sum $t$-intersecting families with $t \ge 2$. If $s$ is even, then again one can consider a $t$-subset $T\subset [n]$ and the corresponding family $$\cF_{n,q,s,r,T}:=\left\{\x\in Q^n:  ~ x_i\ge \frac{s}{2} \mathrm{~for~all} \ i \in T, ~ |S_\x|=r \right\}.$$
In Section 2, we will prove that for large enough $n$ these families contain the largest number of vectors among all $s$-sum $t$-intersecting families.
We will also determine the extremal families if $s$ is odd, but as in that case their definition is more technical, we postpone their introduction to Section 2.


\subsection{Results on intersecting  pairs of vectors
}

We continue with stating Bollob\'as's classical theorem on intersecting set-pair systems for which we prove sum-intersecting analogs. To do so, we recall that if $\cS=\{(A_i,B_i):i=1,2,\dots,n\}$ with $A_i\cap B_i=\emptyset$ for all $1 \le i \le n$, then
\begin{itemize}
     \item $\cS$ is called a \textit{strong ISP-system} (shorthand for intersecting set-pair system) if
    $A_i\cap B_j\neq\emptyset$ for all $1 \le i\neq j \le n$;
     \item $\cS$ is called a \textit{weak ISP-system} if at least one of $A_i\cap B_j\neq\emptyset$ and $B_i\cap A_j\neq\emptyset$ holds for all $1 \le i\neq j \le n$.
 \end{itemize}
If also $a= \max_{1 \le i \le n} |A_i|$ and $b=\max_{1 \le i \le n} |B_i|$, then $\cS$ is a strong or weak $(a,b)$-system.

\begin{theorem}[Bollob\'as \cite{B}]\label{boll}
If $\cS=\{(A_j,B_j):j=1,2,\dots,m\}$ is a strong ISP-system, then the inequality $$\sum_{j=1}^m\frac{1}{\binom{|A_j|+|B_j|}{|A_j|}}\le 1$$ holds. In particular, if $\cS$ is a strong $(a,b)$-system, then $m\le \binom{a+b}{a}$.
\end{theorem}

The following general inequality is valid for weak ISP-systems.

\begin{theorem}[Tuza \cite{T87}]\label{zs87}
Let $0<p<1$ be any real number and $q=1-p$.
If $\{(A_j,B_j):j=1,2,\dots,m\}$ is a weak ISP-system, then the inequality $$\sum_{j=1}^m p^{|A_j|}\, q^{|B_j|}\le 1$$ holds. Moreover, for every $a,b\in\mathbb{N}$ there exists a weak $(a,b)$-system for which equality holds for all $0<p<1$ and $q=1-p$.
\end{theorem}


For a general overview on ISP-systems and their applications in extremal combinatorics we refer to the two-part survey \cite{T1,T2}.
Theorem \ref{zs87} implies the upper bound $m\le \frac{(a+b)^{a+b}}{a^a \, b^b}$ for weak $(a,b)$-systems. 
The best lower bounds on the maximum size of weak $(a,b)$-systems are due to Kir\'aly, Nagy, P\'alv\"olgyi and Visontai \cite{KNPV}, and Wagner \cite{W}.

Now we would like to generalize these notions to vector pairs in the $s$-sum intersecting setting. 
Note that there is no assumption on the size of the ground set of ISP-systems. 
Let us denote by $Q^{<\mathbb{N}}$ $(\subset \mathbb{Z}^{<\mathbb{N}}) $ the set of all $\x \in Q^{\mathbb{N}}$ with finite support $S_\x=\{i:x_i>0\}$.








Assume that for $\cF=\{(\x^j,\y^j) \in Q^{<\mathbb{N}} \times Q^{<\mathbb{N}} : j=1,2,\dots,m \}$ we have $|\x^j\cap_s \y^j|=0$ for all $j=1,2,\dots,m$. We say that is $\cF$ a \textit{strong $s$-sum IVP-system} in $Q^{<\mathbb{N}}$, if  $|\x^i\cap_s \y^j|\neq 0$ for all $1\le i\neq j \le m$ and $\cF$ is a \textit{weak $s$-sum IVP-system} in $Q^{<\mathbb{N}}$, if for all $1\le i\neq j \le m$ \textit{at least one} pair of $\x^i,\y^j$ or $\x^j,\y^i$ is $s$-sum intersecting.
If the supports of all $\x^j$ have size at most $a$, and the supports of all $\y^j$ have size at most $b$, then we will talk about \textit{strong} and \textit{weak $s$-sum $(a,b)$-systems}.









\vskip 0.2truecm

The next observation shows that it is enough to deal with $(q+1)$-sum IVP-systems in $Q^{<\mathbb{N}}$.

\begin{obs}\label{sq}

\  

(i) If $\cF$ is a strong/weak $s$-sum $(a,b)$-system, then for all $(\x,\y)\in \cF$ and all $i\le m$ we have $x_i,y_i<s$. 

(ii) If $\cF\subset Q^{<\mathbb{N}}\times Q^{<\mathbb{N}}$ is a strong/weak $(q+t)$-sum $(a,b)$-system with $t>1$, then there exists a $(q-t+2)$-sum strong/weak $(a,b)$-system $\cF'\subset (\{0,1,\dots,q-t+1\}^{<\mathbb{N}})^2$ with $|\cF|=|\cF'|$.
\end{obs}

\begin{proof}
If $x_i\ge s$ or $y_j\ge s$ for some $(\x,\y)\in \cF$, then $|\x\cap_s \y|>0$. This implies (i).

To see (ii), for any $(\x,\y)\in \cF$ introduce $(\x',\y')$ with $x'_i=\max\{x_i-t+1,0\}$, $y'_i=\max\{y_i-t+1,0\}$ for all indices $i$. 
Clearly, for any $(\x,\y)\in \cF$ and index $j$, we have $x'_j+y'_j<q+t-2(t-1)=q-t+2$. Furthermore, if $|\x^{h_1}\cap_{q+t} \y^{h_2}|>0$, then there exists an index $j$ with $q+t\le x^{h_1}_j+y^{h_2}_j$. So $x^{{h_1}'}_j+y^{{h_2}'}_j \ge q+t-2(t-1)=q-t+2$, and thus the system $\cF'=\{(\x',\y'):(\x,\y)\in \cF\}\subset (\{0,1,\dots,q-t+1\}^{<\mathbb{N}})^2$ is a $(q-t+2)$-sum strong/weak $(a,b)$-system.
\end{proof}

To obtain bounds on the size of $(q+1)$-sum IVP-systems, we write
$m(q,k)$ and $m'(q,k)$ for the maximum number of vector pairs in a strong/weak $(q+1)$-sum $(k,k)$-system.
In particular, for $q=2$ and $s=3$ let $m(k) := m(2,k)$.

To estimate $m(k)$, we let
 $$
   f(k) := \max \frac{(x+y+z)!}{x! \, y! \, z!},
 $$
where the maximum is taken over all nonnegative integers $x,y,z$ such that $x+z\leq k$ and $y+z\leq k$. The following inequalities provide an almost tight bound on $m(k)$,
with only a linear multiplicative error in $k$, while the
function is exponential.

\begin{theorem}\label{q2s3isps-uj}
For every $k\geq 1$ we have
 $$
   f(k) \leq m(k) \leq k \cdot f(k) \,.
 $$
\end{theorem}

Finally, we determine the order of magnitude of the maximum size of strong and weak $(q+1)$-sum IVP systems in $Q^{<\mathbb{N}}$ up to a polynomial factor.

\begin{theorem}\label{weaklimit-uj}
For any $q\ge 1$, $\lim_{k\rightarrow \infty}\sqrt[k]{m(q,k)}=\lim_{k\rightarrow \infty}\sqrt[k]{m'(q,k)}=(\sqrt{q}+1)^2$.
\end{theorem}


A standard calculation shows that the maximum in the definition of $f(k)$ is attained when $z=(1-\frac{1}{\sqrt{2}})k+O(1)$ and $x=y=k-z$. Plugging in these values, we obtain that $f(k)=(c+o(1))\frac{1}{k}(3+2\sqrt{2})^k$ for some real $c<1$. The upper bound of Theorem \ref{q2s3isps-uj} on strong $3$-sum $(k,k)$-systems is a constant factor smaller than the upper bound of Theorem \ref{weaklimit-uj} on weak $3$-sum $(k,k)$-systems. 

\section{Sum-intersecting families of vectors}

This section contains the proofs of Theorem \ref{qekr}, Theorem \ref{rankqekr} and Theorem \ref{qtekr}, but we consider first non-uniform $(q+1)$-sum intersecting vector families.

\begin{proposition}\label{nonuni}
 For $n,q \ge 1$ if $\cF\subseteq Q^n$ is $(q+1)$-sum intersecting, then $|\cF|\le \lceil\frac{(q+1)^n}{2}\rceil$ 
and this bound is best possible.
\end{proposition}

\begin{proof}
Note that we cannot have $\x$ and $\overline{\x}$ both belong to $\cF$. Moreover, there exists one vector $\x$ with $\overline{\x}=\x$ if and only if $q$ is even. This proves the upper bound. For the lower bound consider the family of all vectors with rank larger than $\frac{qn}{2}$ together with one vector from each pair of (the not necessarily different vectors) $\x,\overline{\x}$ of rank $\frac{qn}{2}$ (if such pairs exist).
\end{proof}

\begin{corollary}
For $n,q, s \ge 1$ with $q \ge s$ if $\cF\subseteq Q^n$ is $s$-sum intersecting, then $|\cF|\le (q+1)^n-s^n+ \lceil \frac{s^n}{2} \rceil$ and this bound is best possible.
\end{corollary}

\begin{proof}
If a vector contains an entry at least $s$, then it $s$-sum intersects every other vector. The number of such vectors is $(q+1)^n-s^n$, and then we apply Proposition \ref{nonuni} to the set of all other vectors.
\end{proof}


Now we turn our attention to (rank- or support-) uniform families of vectors. We start 
with the proof of Theorem~\ref{rankqekr}, but we need several definitions and some results from the literature.

\begin{defn}\label{shadows}
The \textit{shadow $\Delta(F)$ of a set $F$} is $\{G\subset F: |G|=|F|-1\}$ and the \textit{shadow $\Delta(\cF)$ of a family $\cF$ of sets} is $\cup_{F\in \cF}\Delta(F)$. If $\cF$ is $r$-uniform and $0 \le \ell < r$, then $\Delta_\ell(\cF):=\{G: |G|=\ell ~\text{and}\ \exists F\in \cF ~\text{s.t. } G \subset F\}$.

\end{defn}

We introduce the notation $<_{\textrm{colex}}$ for the \textit{colex ordering} of all finite subsets of the positive integers. In this ordering for two finite sets $A$ and $B$ we have $A<_{\textrm{colex}}B$ if and only if the largest element of the symmetric difference $(A\setminus B)\cup (B\setminus A)$ of $A$ and $B$ belongs to $B$. 

\vspace{2mm}

Kruskal and Katona independently proved the following fundamental theorem.

\begin{theorem}[Kruskal \cite{kr63}, Katona \cite{k68}]\label{shadow}
Let $n,r,m \ge 1$ and $L_{r,m}$ be the initial segment of $\binom{[n]}{r}$ of size $m$ with respect to the colex ordering. For any $\cF\subseteq \binom{[n]}{r}$ of size $m$, we have $|\Delta(\cF)|\ge |\Delta(L_{r,m})|$.
\end{theorem} 

\vspace{2mm} 

We can introduce the notion of shadow for vectors, too.

\begin{defn}\label{vectshadow}
The \textit{shadow $\Delta(\x)$ of a vector $\x \in Q^n$} is $\{\y< \x: r(\y)=r(\x)-1\}$, where $<$ denotes the coordinate-wise ordering, i.e., for two vectors $\x$ and $\y$ we have $\y < \x$ if and only if $y_i \le x_i $ for all $1 \le i \le n$ and $y_i < x_i$ for at least one $i$. Then for $\cF \subseteq Q^n$ we define the \textit{shadow $\Delta(\cF)$ of $\cF$} as $\cup_{\x\in \cF}\Delta(\x)$ and for $r$-rank uniform $\cF$ and $\ell<r$ we let $\Delta_\ell(\cF)=\{\y: r(\y)=\ell ~\text{and}\ \exists \x\in \cF ~\text{such that } \y < \x\}$. We will write $\x\leqslant \y$ for $\x<\y$ or $\x=\y$.
\end{defn}


Analogously to the set case, we can introduce the \textit{colex ordering of $Q^n$}, i.e., for $\x,\y \in Q^n$ we have $\x<_{\textrm{colex}}\y$ if and only if $x_i<y_i$ where $i$ is the largest coordinate in which $\x$ and $\y$ differ. 

\vspace{2mm}

Clements and Lindstr\"om proved a generalization of the Kruskal-Katona theorem for the shadows of vectors introduced in Definition \ref{vectshadow}.

\begin{theorem}[Clements, Lindstr\"om \cite{cl69}]\label{qshadow}
Let $q,r,m,n \ge 1$, and let $L_{q,r,m}$ be the initial segment of $Q(n,r)$ of size $m$ with respect to the colex ordering. For any $\cF\subseteq Q(n,r)$ of size $m$, we have $|\Delta(\cF)|\ge |\Delta(L_{r,m})|$.
\end{theorem} 

One can easily check the following properties of the colex ordering of sets and vectors, so we omit their proof. 

\begin{proposition}\label{colex} Suppose $n \ge r \ge 1$. \

\vspace{1mm}

(i) Both in $\binom{[n]}{r}$ and in $Q(n,r)$, the shadow of an initial segment is an initial segment, so one can iterate Theorems \ref{shadow} and \ref{qshadow} to obtain that initial segments minimize the size of shadows of any lower rank.

(ii) If $\cF$ is the family of the largest $m$ sets of $\binom{[n]}{r}$ with respect to the colex ordering, then $\overline{\cF}=L_{n-r,m}$.

(iii) If $\cF$ is the family of the largest $m$ vectors of $Q(n,r)$ with respect to the colex ordering, then $\overline{\cF}=L_{q,qn-r,m}$.
\end{proposition}

Before the proof of Theorem \ref{rankqekr}, let us briefly recall the proof of the upper bound in Theorem \ref{ekr} that uses the Kruskal--Katona shadow theorem (Theorem \ref{shadow}) and was obtained by Daykin \cite{D} as we would like to mimic it.

Suppose contrary to the statement of Theorem \ref{ekr} that there exists an intersecting family $\cF\subseteq \binom{[n]}{r}$ of size larger than $\binom{n-1}{r-1}$. Consider the family $\overline{\cF}=\{[n]\setminus F:F\in \cF\}$ and observe that as $\cF$ is intersecting and $n \ge 2r$, we must have $\cF\cap \Delta_r(\overline{\cF})=\emptyset$.
Clearly, $|\overline{\cF}|=|\cF|>\binom{n-1}{r-1}=\binom{n-1}{n-r}$. Applying Theorem \ref{shadow}, any $(n-r)$-uniform family of size larger than $\binom{y}{n-r}$ has $r$-shadow larger than $\binom{y}{r}$. So $\binom{n}{r}=|\binom{[n]}{r}|\ge |\cF|+|\Delta_r(\overline{\cF})|>\binom{n-1}{r-1}+\binom{n-1}{r}=\binom{n}{r}$. This contradiction proves the upper bound in Theorem \ref{ekr}. 

This proof seems to be very lucky that it includes miraculous equalities $\binom{n-1}{r-1}=\binom{n-1}{n-r}$ and $\binom{n-1}{r-1}+\binom{n-1}{r}=\binom{n}{r}$, so let us recite it without any calculation. Consider greedily the largest sets of $\binom{[n]}{r}$ with respect to the colex order as long as they form an intersecting family. Let $\cF_0$ be the family when we need to stop. If $\cF_0 \cup \Delta_r(\overline{\cF_0})=\binom{[n]}{r}$, then $\cF_0$ is a largest possible intersecting family. Indeed, if $|\cF|>|\cF_0|$, then as $\overline{\cF_0}$ is an initial segment, by Proposition \ref{colex} (i) and (ii), we have $|\cF|+|\Delta_r(\overline{\cF})|> |\cF_0|+|\Delta_r(\overline{\cF_0})|=\binom{n}{r}$, so $\cF$ cannot be intersecting. To obtain the results of Theorem \ref{ekr} about intersecting families, all we need to observe is that $\cF_0=\{F\in \binom{[n]}{r}:n\in F\}$ and $\Delta_r(\overline{\cF_0})=\{F\in \binom{[n]}{r}:n\notin F\}$. 

\vspace{2mm}

Before the proof of Theorem \ref{rankqekr} let us restate it.

\begin{reptheorem}{rankqekr}
Let $n,q,r \ge 1$ and $\cF\subseteq Q^n$ be an $r$-rank uniform $(q+1)$-sum intersecting family with $\frac{q+1}{2}\le r\le \frac{qn}{2}$. Then
\begin{eqnarray}
|\cF|\le\left\{
\begin{array}{cc} 
\sum_{j=\frac{q+1}{2}}^q|Q(n-1,r-j)| & \textnormal{if} ~q+1 ~\textnormal{is even}, \\
1+\sum_{j=\lceil\frac{q+1}{2}\rceil}^q\sum_{i=1}^{\lfloor \frac{2(r-1)}{q}\rfloor}|Q(n-i,r-j-\frac{(i-1)q}{2})|& \textnormal{if $q+1$ is odd},
\end{array}
\right.
\end{eqnarray}
and these bounds are best possible.

\end{reptheorem}

\begin{proof}
Clearly $\x \in Q^n$ does not $(q+1)$-sum intersect a vector $\y \in Q^n$ if and only if $\y\leqslant \overline{\x}$. Also, $\cF \subseteq Q(n,r)$ is a $(q+1)$-sum intersecting family if and only if $\cF \cap \Delta_r(\overline{\cF})$
contains at most one vector as $\cF$ may contain one vector $\x$ that does not $(q+1)$-sum intersect itself. Indeed, if $\x\neq \y$ and $|\x\cap_{q+1}\y|=0$, then $\x,\y\in \cF \cap \Delta_r(\overline{\cF})$. On the other hand if $\x,\y\in \cF\cap \Delta_r(\overline{\cF})$ and $\cF$ is intersecting, then by the above, we cannot have $\x\leqslant\overline{\y}$. Therefore, we must have $\x\leqslant \overline{\x}$ and $\y\leqslant\overline{\y}$. But as $|\x\cap_{q+1} \y|>0$, there must exist an index $i$ with $x_i+y_i\ge q+1$, so either $x_i$ or $y_i$, say $x_i$, is at least $\frac{q+1}{2}$. But then $x_i>q-x_i=\overline{x}_i$, a contradiction.

The reasoning of Daykin stays valid with a little modification, if for the maximal $(q+1)$-sum intersecting family $\cF_0\subseteq Q(n,r)$ consisting of largest vectors with respect to the colex ordering we have both $\cF_0 \cup \Delta_r(\overline{\cF_0})=Q(n,r)$ \textit{and} $|\Delta_r(\overline{\cF}_0)|<|\Delta_r(\overline{\cF}_0^+)|$, where $\overline{\cF}_0^+$ is the initial segment of the colex ordering of $Q(n,qn-r)$ one larger than $\overline{\cF}_0$. Indeed, if $\cF$ was an $r$-rank uniform $(q+1)$-sum intersecting family larger than $\cF_0$, then we would get a contradiction by the following series of inequalities: $$|\cF\cup \Delta_r(\overline{\cF})|\ge |\cF|+|\Delta_r(\overline{\cF})|-1\ge |\cF_0|+1+|\Delta_r(\overline{\cF}_0)|+1-1=|Q(n,r)|+1.$$ 

\vspace{3mm}

And this is exactly the case: for the maximal $(q+1)$-sum intersecting family $\cF_0\subseteq Q(n,r)$ consisting of largest vectors with respect to the colex ordering, we prove that we have both $\cF_0 \cup \Delta_r(\overline{\cF_0})=Q(n,r)$ \textit{and} $|\Delta(\overline{\cF}_0)|<|\Delta(\overline{\cF}_0^+)|$, where $\overline{\cF_0}^+$ is the one larger initial segment of the colex ordering of $Q(n,qn-r)$ than $\overline{\cF}_0$.

\vspace{2mm}

Suppose first that $q+1=2k$. Then $\cF_0=\{ \x \in Q(n,r):x_n\ge k\}$, $\overline{\cF}_0=\{ \x \in Q(n,qn-r):x_n<k\}$ and clearly $\Delta_r(\overline{\cF}_0)=\{ \x \in Q(n,r):x_n<k\}=Q(n,r)\setminus \cF_0$; and since $\overline{\cF}_0^+$ contains a vector $\x$ with $x_n=k$, its $r$-shadow is strictly larger than that of $\overline{\cF}_0$.

\vspace{2mm}

Suppose next $q+1=2k+1$. Then $$\cF_0=\bigcup_{j=0}^{\lfloor \frac{r-1}{k}\rfloor}\{ \x \in Q(n,r):x_n=x_{n-1}=\dots =x_{n-j+1}=k, \ x_{n-j}>k\} \cup \{ \x^*\},$$
where $x^*_n=x^*_{n-1}=\dots =x^*_{n-\lfloor \frac{r-1}{k}\rfloor}=k, \ x^*_{n-\lfloor \frac{r-1}{k}\rfloor-1}\equiv r ~ (\mathrm{mod} ~k)$ and all other entries are 0. Observe that $\x^*$ does not $(q+1)$-sum intersect itself. To see that $\cF_0 \cup \Delta_r(\overline{\cF}_0)=Q(n,r)$ holds, one only has to observe that any vector $\y \in Q(n,r)\setminus \cF_0$ with $y_n=y_{n-1}=\dots =y_{n-\lfloor \frac{r-1}{k}\rfloor}=k$ belongs to $\Delta_r(\overline{\x^*})$. Also, any vector $\y \in Q(n,r)$ with $\overline{\x^*}<_{\textrm{colex}} \y$ has an entry larger than $k$ in the last $\lfloor\frac{r-1}{k}\rfloor$ coordinates, so $|\Delta_r(\overline{\cF}_0^+)|>|\Delta_r(\overline{\cF}_0)|$.

\vspace{2mm}

This completes the proof of Theorem \ref{rankqekr}.
\end{proof}

\vskip 0.2truecm

We continue with the proof of Theorem \ref{qekr}. Before doing so, we cite two well-known stability results that we use during the proof of Theorems \ref{qekr} and \ref{qtekr}.

\begin{theorem}[Hilton, Milner \cite{HM}]\label{hm}
If $\cF\subseteq \binom{[n]}{r}$ is an intersecting family with $n\ge 2r+1$ and $\cap_{F\in\cF}F=\emptyset$, then $|\cF|\le \binom{n-1}{r-1}-\binom{n-r-1}{r-1}+1$.
\end{theorem}

\begin{theorem}[Frankl \cite{F}]\label{frankl}
Let $\cF\subseteq \binom{[n]}{r}$ be a $t$-intersecting family with $|\cap_{F\in \cF}F|<t$. If $n$ is large enough, then $|\cF|\le \max\{|\cF_1|,|\cF_2|\}$, where
$$\cF_1=\left\{F\in\binom{[n]}{r}:[t]\subset F, F\cap [t+1,r+1]\neq \emptyset \right\}\cup \binom{[r+1]}{r}$$
and
$$\cF_2=\left\{F\in\binom{[n]}{r}:|F\cap [t+2]|\ge t+1\right\}.$$
\end{theorem}

Now let us restate Theorem \ref{qekr}.

\begin{reptheorem}{qekr}

For any $s>q\ge 2$ and integer $r \ge 1$,
if $\cF\subseteq Q^n$ is $r$-support uniform $s$-sum intersecting with $n\ge qr^2$, then 
\begin{eqnarray}
|\cF|\le\left\{
\begin{array}{cc} 
(q-\frac{s}{2}+1) q^{r-1}\binom{n-1}{r-1} & \textnormal{if} ~s ~\textnormal{is even,} \\
1+(q-\lceil\frac{s}{2}\rceil+1)\sum_{i=1}^r \binom{n-i}{r-i}q^{r-i}& \textnormal{if $s$ is odd,}
\end{array}
\right.
\end{eqnarray}
and these bounds are best possible.

\end{reptheorem}

\begin{proof}[Proof of Theorem \ref{qekr}]
Suppose first  that $s$ is even. The constructions  showing that the bound is best possible are $\cF_{n,q,s,r,i}=\{ \x \in Q^n: \frac{s}{2}\le x_i\}$. To see the upper bound, let $\cF$ be an $r$-support uniform $s$-sum intersecting family and let $\cS_\cF$ denote the family of supports in $\cF$. For a fixed support $S$, the number of vectors having $S$ as support is bounded by a constant (depending on $r$ and $q$), therefore, by Theorem \ref{hm}, unless all supports in $\cS_\cF$ share a common element $i$, we have $|\cF|\le q^rr\binom{n-2}{r-2}<(q-\frac{s}{2}+1)q^{r-1}\binom{n-1}{r-1}$ if $n\ge qr^2$. So we can suppose that there exists  an index $i$ that belongs to all supports. Assume next that there exists $\x \in \cF$ with $x_i<\frac{s}{2}$. Then consider the subfamily $\cF'=\{ \y \in \cF: y_i\le \frac{s}{2}\}$. As vectors in $\cF'$ must all $s$-sum intersect $\x$, but they do not $s$-sum intersect it at coordinate $i$, therefore their supports must intersect the support of $\x$ in some coordinate other than $i$. Therefore, we obtain $|\cS_{\cF'}|\le (r-1)\binom{n-2}{r-2}$ and thus $|\cF'|\le \frac{s}{2}q^{r-1}(r-1)\binom{n-2}{r-2}$. But then 
$$|\cF|\le |\cF'|+(q-\frac{s}{2})q^{r-1}\binom{n-1}{r-1}<(q-\frac{s}{2}+1)q^{r-1}\binom{n-1}{r-1}$$ 
if $n\ge r^2\frac{s}{2}$. We obtained that either $\cF$ is smaller than the claimed bound or $\cF\subseteq \cF_{n,q,s,r,i}$ for some index $i$.

Suppose next that $s$ is odd. The extremal families are defined via ordered $r$-tuples $(i_1,i_2, \dots, i_r)$ in the following way:
$$ \cF_{n,q,s,(i_1,i_2,\dots,i_r)}=\{ \x \}\cup \bigcup_{j=1}^r\{\y \in Q^n: y_{i_1}=y_{i_2}\dots =y_{i_{j-1}}=\lfloor \frac{s}{2}\rfloor, \ y_{i_j}\ge \frac{s}{2}\},$$ 
where $\x$ is the vector with $x_{i_j}=\lfloor \frac{s}{2}\rfloor$ for all $1\le j\le r$ and $x_i=0$ otherwise. To prove the upper bound, we proceed by induction on $r$. If $r=1$, then all supports of an $r$-support uniform $s$-sum intersecting family $\cF$ must be the same singleton $\{i\}$. If $m$ is the minimum entry over all vectors in $\cF$ at coordinate $i$, then all other entries must be at least $s-m$, so the number of vectors is at most $\min\{q-m+1,q-(s-m)\}$. This is maximized if $m=\lfloor \frac{s}{2}\rfloor$ and the claimed bound follows. Let $r>1$, and $\cF\subseteq Q^n$ be an $r$-support uniform, $s$-sum intersecting family. Then just as in the even $s$ case, using Theorem \ref{hm}, we obtain that $|\cF|\le q^rr\binom{n-2}{r-2}<q^{r-1}\binom{n-1}{r-1}$ unless all sets in $\cS_\cF$ share a common element $i_1$ or $n\le qr^2$. If there exists a vector $\z \in \cF$ with $z_{i_1}<\lfloor \frac{s}{2}\rfloor$, then also just as in the even $s$ case, we obtain that $\cF'=\{ \y \in \cF: y_{i_1}\le \lceil \frac{s}{2} \rceil\}$ is of size at most $\lceil\frac{s}{2}\rceil q^{r-1}(r-1)\binom{n-2}{r-2}$ and thus $\cF$ is smaller than the claimed bound if $n\ge sr^2$. So we can assume that for all vectors $\z \in \cF$, we have $z_{i_1}\ge \lfloor \frac{s}{2}\rfloor$. The number of those vectors $\z$ with $z_{i_1}\ge \lceil \frac{s}{2} \rceil$ is $(q-\lceil \frac{s}{2}\rceil +1)q^{r-1}\binom{n-1}{r-1}$, while the family $\cF^*=\{\z'\in \cF:z_{i_1}=\lfloor \frac{s}{2}\rfloor\}$ is $(r-1)$-support uniform, $s$-sum intersecting, where $\z'$ is the vector obtained from $\z$ by removing its $i_1$st entry. As $n-1 \ge qr^2-1\ge q(r-1)^2$, by induction, we obtain $$|\cF^*|\le 1+(q-\lceil\frac{s}{2} \rceil+1)\sum_{i=1}^{r-1}q^{r-1-i}\binom{n-1-i}{r-1-i}$$ and so 

\begin{eqnarray*}
    |\cF| & \le & |\cF^*|+(q-\lceil \frac{s}{2}\rceil +1)q^{r-1}\binom{n-1}{r-1} \\
    & \le & 1+(q-\lceil\frac{s}{2}\rceil+1)\sum_{i=1}^{r-1}q^{r-1-i}\binom{n-1-i}{r-1-i}+(q-\lceil \frac{s}{2} \rceil +1)q^{r-1}\binom{n-1}{r-1} \\
    & = & 1+(q-\lceil\frac{s}{2}\rceil+1)\sum_{i=1}^{r}q^{r-i}\binom{n-i}{r-i},
\end{eqnarray*}
as claimed.
\end{proof}

In the remainder of this section, we consider $s$-sum $t$-intersecting families with $t\ge 2$.

\begin{con}\label{cont}
For any $n,q,r,t \ge 1$ with $n \ge r \ge t$ and $s$ even with $q< s \le 2q$ and for any $T \in \binom{[n]}{t}$ let us define
$$\cF_{n,q,s,r,T}:=\left\{\x\in Q^n:  ~ x_i\ge \frac{s}{2} \mathrm{~for~all} \ i \in T, ~ |S_\x|=r \right\}.$$
Observe that the size of $\cF_{n,q,s,r,T}$ is $(q-\frac{s}{2}+1)^t q^{r-t}\binom{n-t}{r-t}$.

\smallskip

For $n,r,q,t \ge 1$ with $n \ge r\ge t$ and $s$ odd with $q<s<2q$ let us define the following $r$-support uniform families: for any $T'\in  \binom{[r]}{t-1}$ we pick a vector $\x_{T'}$ with $S_{\x_{T'}}=[r]$ such that $(x_{T'})_i>\frac{s}{2}$ for all $i\in T'$, and $(x_{T'})_i=\lfloor \frac{s}{2}\rfloor$ for all $i\in [r]\setminus T'$. Then we have

    $$\cF_{n,q,s,r,t}  := \left\{\x_{T'}: T'\in \binom{[r]}{t-1}\right\}\cup$$
    $$\bigcup_{T \in \binom{[r]}{t}}\left\{\y \in Q^n: (\forall i\in T)(y_{i}>\lfloor \frac{s}{2}\rfloor) \wedge (\forall i \in [\max T]\setminus T) (y_{i}=\lfloor\frac{s}{2}\rfloor) \wedge |S_{\y}|=r \right\}.$$
In words, the vectors belonging to the second row of the definition have at least $t$ indices $i_1<i_2<\dots <i_t\le r$ such that their corresponding entries have value strictly greater than $\lfloor s/2\rfloor$, and all entries with indices smaller than $i_t$ have value at least $\lfloor s/2\rfloor$.

\smallskip

The family $\cF_{n,q,s,r,t}$ is $s$-sum $t$-intersecting as for any pair $\x,\y \in \cF_{n,q,s,r,t}$ there exist at least $t$ coordinates $i \in [r]$, where one of $x_i,y_i$ is at least $\lfloor \frac{s}{2}\rfloor$ while the other is at least $\lceil \frac{s}{2}\rceil$.
\end{con}

Let $f(n,q,s,r,t)$ denote the size of $\cF_{n,q,s,r,t}$.

\begin{theorem}\label{qtekr}
For any $2\le q<s\le 2q$ and $r \ge t \ge 1$, there exists $n(q,s,r,t)$ such that if $\cF\subseteq Q^n$ is $r$-support uniform $s$-sum $t$-intersecting with $n\ge n(q,s,r,t)$, then 
\begin{eqnarray}
|\cF|\le\left\{
\begin{array}{cc} 
(q-\frac{s}{2}+1)^t q^{r-t}\binom{n-t}{r-t} & \textnormal{if} ~s ~\textnormal{is even}, \\
f(n,q,s,r,t) & \textnormal{if $s$ is odd},
\end{array}
\right.
\end{eqnarray}
and these bounds are best possible as shown by the families of Construction \ref{cont}.
\end{theorem}

We will need the following simple observations on $f(n,q,s,r,t)$.

\begin{proposition}\label{contprop} Suppose that $n,q,s,r,t$ are integers with the assumptions on them as in Construction~\ref{cont}.

\vspace{2mm}

(i) If $r\ge 2t$, then $$f(n,q,s,r,t)=$$ $$\binom{n-t}{r-t}q^{r-t}(q-\lfloor \frac{s}{2}\rfloor)^t+\sum_{S\subsetneq [t]}(q-\lfloor \frac{s}{2}\rfloor)^{|S|}f(n-t,q,s,r-t,t-|S|).$$

(ii) If $t<r<2t$, then $$f(n,q,s,r,t)=$$ $$\binom{n-t}{r-t}q^{r-t}(q-\lfloor \frac{s}{2}\rfloor)^t+\binom{t}{2t-r-1}+\sum_{S\subsetneq [t], |S|\ge 2t-r}(q-\lfloor \frac{s}{2}\rfloor)^{|S|}f(n-t,q,s,r-t,t-|S|).$$

(iii) If $r>t$, then $f(n,q,s,r,t)>(q-\lfloor s/2\rfloor)^tq^{r-t}\binom{n-t}{r-t}$.
\end{proposition}

\begin{proof}
In all of (i), (ii) and (iii), the first term of the right-hand side (in case of (iii), the only term) stands for those vectors for which $x_i >\frac{s}{2}$ for all $1\le i \le t$. In (i), the big sum partitions the other vectors according to which of the first $t$ entries have value greater than $s/2$. 

In (ii), the big summation can neglect small subsets $S$ of $[t]$, because if $|S|<2t-r$, then for any $\x\in \cF_{n,q,s,r,t}$ with $\{i\in [t]: x_i>s/2\}=S$, the number of indices $i$ in $[r]$ for which $x_i>s/2$ is at most $|S|+r-t<t$. So if $|S|<2t-r$, then to reach at least $t-1$ such indices (the minimum for a vector in $\cF_{n,q,s,r,t}$), we need exactly $2t-r-1$ indices from $[t]$. The middle term stands for those $\x_{T'}$s where $T'$ contains exactly $2t-r-1$ elements from $[t]$.
\end{proof}

Now we continue with the proof of Theorem \ref{qtekr}.

\begin{proof}[Proof of Theorem \ref{qtekr}]

Suppose first that $s$ is even. To see the upper bound, let $\cF$ be an $r$-support uniform $s$-sum $t$-intersecting family and let $\cS_\cF$ denote the family of supports in $\cF$. For a fixed support $S$, the number of vectors having $S$ as support is bounded by a constant (depending on $r$ and $q$), therefore, by Theorem \ref{frankl}, unless all supports in $\cS_\cF$ share all
elements of a $t$-subset $T$ of $[n]$, we have $|\cF|=O(n^{r-t-1})<\binom{n-t}{r-t}$ if $n$ is large enough. So we can suppose that there exists  a $t$-subset $T$ that is contained in all supports. Assume next that there exists $\x \in \cF$ with $x_i<\frac{s}{2}$ for some $i \in T$. Then consider the subfamily $\cF'=\{ \y \in \cF: y_i\le \frac{s}{2}\}$. As vectors in $\cF'$ must all $s$-sum $t$-intersect $\x$, but they do not $s$-sum intersect it at coordinate $i$, therefore their supports must intersect the support of $\x$ in some coordinate outside $T$. Therefore, we obtain $|\cF'|=O(n^{r-t-1})$. But then $$|\cF|\le |\cF'|+(q-\frac{s}{2})(q-\frac{s}{2}+1)^{t-1}q^{r-t}\binom{n-t}{r-t}<(q-\frac{s}{2}+1)^t q^{r-t} \binom{n-t}{r-t}$$ if $n$ is large enough. We obtained that either $\cF$ is smaller than the claimed bound or $\cF\subseteq \cF_{n,q,s,r,T}$ for some $t$-subset $T$.

\vspace{3mm}

Suppose next that $s$ is odd. We proceed by induction on $r+t$ and observe that in all cases, the family of supports must be $t$-intersecting. The case $t=1$ is covered by Theorem \ref{qekr}. Let $\cF\subseteq Q^n$ be an $s$-sum $t$-intersecting $r$-support uniform family. We consider three cases according to the relationship of $r$ and $t$.

\vskip 0.2truecm

\textsc{Case I:} $r=t.$

\vskip 0.2truecm

The assumption $r=t$ implies that all supports in $\cF$ are identical, say the support is $S$. Therefore, for any $\x, \y \in \cF$ and $i \in S$ we must have $x_i+y_i\ge s$. In particular, for any $i \in S$ there is at most one $\x \in \cF$ with $x_i<s/2$. So $$|\cF|\le \binom{t}{t-1}+(q-\lfloor s/2\rfloor)^{t},$$ as claimed.

\vskip 0.2truecm

\textsc{Case II:} $t<r$

\vskip 0.2truecm

The family $\cS_\cF$ of supports is $t$-intersecting, so unless all supports of $\cF$ share $t$ elements, we have $|\cF|=O(n^{r-t-1})$ by Theorem \ref{frankl}. Let $T$ be the set of these $t$ elements, and for any $S\subset T$ let $\cF_S$ denote the family of those vectors $\x \in \cF$ for which $x_i\ge s/2$ for all $i\in S$, and $1 \le x_i \le s/2$ for all $i\in T\setminus S$. As all supports contain $T$, we have $\cF=\cup_{S\subset T}\cF_S$. Clearly, $|\cF_T|\le q^{r-t}\binom{n-t}{r-t}(q-\lfloor s/2\rfloor)^t$.

We claim that if there exists $\x\in \cF$ with $x_i<\lfloor s/2\rfloor$ for some $i\in T$, then $|\cF|<f(n,q,s,r,t)$. Indeed, the vectors $\y\in \cF_T$ with $y_i=\lceil s/2\rceil$ $s$-sum $t$-intersect $\x$, so $S_\x\cap S_\y\cap ([n]\setminus T)\neq \emptyset$. Therefore the number of such vectors is $O(n^{r-t-1})$ and thus $|\cF_T|\le (q-\lceil s/2\rceil)(q-\lfloor s/2\rfloor)^{t-1}q^{r-t}\binom{n-t}{r-t}+O(n^{r-t-1})$. 
Also, for any $j\in T$, the number of vectors $\z\in \cF$ with $z_j\le \lfloor s/2\rfloor$ is at most $O(n^{r-t-1})$ as to have $|\z\cap_s\z'|\ge t$ for two such vectors, $S_\z$ and $S_{\z'}$ must intersect outside $T$. 
Adding up for all $j\in T$, we obtain $$|\cF|\le (q-\lceil s/2\rceil)(q-\lfloor s/2\rfloor)^{t-1}q^{r-t}\binom{n-t}{r-t}+O(n^{r-t-1})+O(tn^{r-t-1})<f(n,q,s,r,t)$$ as claimed. Here the last inequality follows from Proposition \ref{contprop} (iii) (watch out for floor and ceiling signs!). So we can assume that $x_i\ge \lfloor s/2\rfloor$ for all $i\in T$, thus $x_i=\lfloor s/2\rfloor$ for all $i\in T\setminus S$ and $\x\in \cF_S$. This implies $|\cF_S|\le (q-\lfloor s/2\rfloor)^{|S|}|\cF'_S|$ with $\cF'_S=\{\x': \x \in \cF_S\}$, where $\x'$ is the vector obtained from $\x$ by deleting the coordinates belonging to $T$.

\vskip 0.2truecm

\textsc{Case IIa:} $t<r<2t.$

\vskip 0.2truecm

Consider families $\cF_S$ for all subsets $S$ with $2t-r\le |S|<t$. Observe $\cF'_S$ is $(r-t)$-support uniform $s$-sum $(t-|S|)$-intersecting, and thus by induction, we have $$|\cF_S|\le (q-\lfloor s/2\rfloor)^{|S|}|\cF'_S|\le (q-\lfloor s/2\rfloor)^{|S|}f(n-t,q,s,r-t,t-|S|).$$

Finally, consider all subsets $S\subset T$ with $|S|<2t-r$. As for two vectors $\x',\x''\in \cF_S$, we have $|\x'\cap_s\x''|\le |S|+r-t < 2t-r+r-t=t$, we must have $|\cF_S|\le 1$ for all such $S$. Observe that for any $(r+1-t)$-subset $Z\subset T$ there exists at most one subset $S\subset T$ with $Z\cap S=\emptyset$ and $\cF_S\neq \emptyset$. Indeed, if $\x \in \cF_S$, $\y \in \cF_{S'}$, then $\x$ and $\y$ can only $s$-sum intersect in at most $r-t$ coordinates outside $T$ and in at most $t-(r+1-t)=2t-r-1$ coordinates within $T$, so $|\x \cap_s \y|\le t-1$, a contradiction. Therefore $$\sum_{S\subset T, |S|<2t-r}|\cF_S|\le \binom{t}{r+1-t}=\binom{t}{2t-r-1}.$$ 
Adding up these bounds  for all $|\cF_S|$ together with the bound on $|\cF_T|$, we obtain the desired bound on $|\cF|$ by Proposition~\ref{contprop}~(ii).

\vskip 0.2truecm

\textsc{Case IIb:} $2t\le r.$

\vskip 0.2truecm
 
 In this case, for any $S\subsetneq T$, the family $\cF'_S$ is $(r-t)$-support uniform $s$-sum $(t-|S|)$-intersecting, and thus by induction, we have $$|\cF_S|\le (q-\lfloor s/2\rfloor)^{|S|}|\cF'_S|\le (q-\lfloor s/2\rfloor)^{|S|}f(n-t,q,s,r-t,t-|S|).$$ Adding up these bounds  for all $|\cF_S|$ together with the bound on $|\cF_T|$, we obtain the desired bound on $|\cF|$ by Proposition~\ref{contprop}~(i). 
\end{proof}

We did not elaborate on the value of the threshold $n(q,s,r,t)$. The statement of Theorem \ref{frankl} was proved by Ahlswede and Khachatrian \cite{ak2} under the weaker condition $n\ge (t+1)(k-t+1)+1$. Computations similar to the one in the proof of Theorem \ref{qekr}, one would obtain that the choice $n(q,s,r,t)=q^tr(r+t)$  works.

\section{Intersecting vector pairs}

\def \scap {\cap_{\mathrm{sum}}}

In this section we provide proofs for Theorem \ref{q2s3isps-uj} and Theorem \ref{weaklimit-uj}.


\medskip 

Let us start with a general construction.

\begin{con}
   \label{kon:abc}
Let $c\le a\le b$ and $3 \le s < 2q$ be integers and fix a set $X$ of size $a+b-c$. For any 3-partition $A\cup B\cup C=X$ with $|A|=a-c$, $|B|=b-c$, $|C|=c$, we define the pairs $\x^{A,B,C}$ and $\y^{A,B,C}$ with \[x^{A,B,C}_i=y^{A,B,C}_i=\lceil s/2\rceil-1 ~\text{if}\ i \in C,
\]
\[
x^{A,B,C}_i=\lfloor s/2\rfloor+1, ~ y^{A,B,C}_i=0 ~\text{if}\ i\in A
\]and
\[x^{A,B,C}_i=0, ~ y^{A,B,C}_i=\lfloor s/2\rfloor+1 ~\text{if}\ i\in B.
\]
Note that $\{(\x^{A,B,C},\y^{A,B,C}): A\cup B\cup C=X, \ |A|=a-c, \ |B|=b-c, \ |C|=c\}$ is a strong $s$-sum IVP-system of cardinality  $\binom{a+b-c}{b-c}\binom{a}{c}.$ 

More generally, let $\alpha_0,\alpha_1,\dots,\alpha_q$ be positive integers with $\sum_{i=0}^{q-1}\alpha_i\le a$ and $\sum_{i=1}^{q}\alpha_i\le b$. Set $N=\sum_{i=0}^{q}\alpha_i\le a$, and define
$$
\{(\x^{A_0,A_1,\dots,A_q},\y^{A_0,A_1,\dots,A_q}):[N]=\bigsqcup_{i=0}^qA_i, ~|A_i|=\alpha_i\} \ ,
$$
where $x^{A_0,A_1,\dots,A_q}_j=q-y^{A_0,A_1,\dots,A_q}_j=i$ if and only if $j\in A_i$.

Observe that the above is a strong $(a,b)$-system. Indeed, by definition
we have that $x^{A_0,A_1,\dots,A_q}_j+y^{A_0,A_1,\dots,A_q}_j=q$ for any $j\in N$, and $A_0,A_1,\dots,A_q$ partition $N$ with $|A_i|=\alpha_i$, and so $|\x^{A_0,A_1,\dots,A_q}\cap_{q+1} \y^{A_0,A_1,\dots,A_q}|=0$.  
Furthermore, if $(A_0,A_1,\dots,A_q)\neq (B_0,B_1,\dots,B_q)$, then there exists $j$ such that $A_j\neq B_j$. We consider such $j$ that minimizes $\min\{j,q-j\}$ and we can suppose without loss of generality that $j \le q/2$. By the assumption on $j$, 
there exist $i\in B_j\setminus A_j$ and $i'\in A_j\setminus B_j$. Then we have $x^{A_0,A_1,\dots,A_q}_i+y^{B_0,B_1,\dots,B_q}_i> j + q-j$, as $y^{B_0,B_1,\dots,B_q}_i=j$, $i\in B_j\setminus A_j$ and $j$ is minimal; and we also have $x^{B_0,B_1,\dots,B_q}_{i'}+y^{A_0,A_1,\dots,A_q}_{i'}>q-j+j$ by similar reasons. This proves that we indeed defined a strong $(a,b)$-system.
\end{con}











\subsection{Upper bound for strong 3-sum IVP-systems in $\{0,1,2\}^{<\mathbb{N}}$}

In this subsection we will prove Theorem \ref{q2s3isps-uj}. Let $\{(\x^j,\y^j) \mid 1\leq j\leq m\}$ be a strong 3-sum $(k,k)$-system in $\{0,1,2\}^{<\mathbb{N}}$. Let us also introduce the following further notation for $j=1,\dots,m$: 

\begin{itemize}
    \item $a_j = |S_{\x^j} \setminus S_{\y^j}|$,
    \item $b_j = |S_{\y^j} \setminus S_{\x^j}|$,
    \item $c_j = |S_{\x^j} \cap S_{\y^j}|$.
\end{itemize}

\noindent 
First we prove the following LYM-type theorem for 3-sum $(a,b)$-systems. 

\begin{theorem}
  \label{t:sum-abc}
Suppose that $a,b,m \ge 1$ and $\{(\x^j,\y^j) \mid 1\leq j\leq m\}$ is a strong 3-sum $(a,b)$-system in $\{0,1,2\}^{<\mathbb{N}}$. Then
\begin{equation}\label{eq:sum-abc}
  \sum_{j=1}^m \frac{a_j! \, b_j! \, c_j!}{(a_j+b_j+c_j)!}=\sum_{j=1}^m \frac{1}{\binom{a_j+b_j+c_j}{a_j+b_j} \binom{a_j+b_j}{a_j}} \leq
\min(a,b).
  \end{equation}
\end{theorem}

\begin{proof}
Essentially we apply induction on $n$, that is the size of union of the supports of elements in $\{(\x^j,\y^j) \mid 1\leq j\leq m\}$.
\begin{enumerate}

\medskip 
    
    \item Note first that $a_i=0$ and $b_j=0$ cannot hold simultaneously for any $1 \le i\neq j \le m$.
    Indeed, if $S_{\x^i}\subset S_{\y^i}$ and $S_{\y^j}\subset S_{\x^j}$ then all nonzero entries in $\x^i$ are equal to 1, and the same holds for all nonzero entries in $\y^j$ as well, hence $\x^i \cap_3 \y^j =\emptyset$, a contradiction.
    As a consequence, either $a_j>0$ for all $j$ or $b_j>0$ for all $j$ (or both), or there is exactly one $j$ with $a_j=b_j=0$.

\medskip 

    \item As long as $S_{\y^j}\not\subseteq S_{\x^j}$ holds for all $j$:

\smallskip 
    
    For every $t \in [n]$, consider the systems
    $$
      \{ (\x^j, (\y^j)') \mid 1\leq j\leq m, \ t \notin S_{\x^j} \},
    $$ 
      
    where $((y^j)')_i = (y^j)_i$ for all $i \in [n]\setminus\{t\}$ and $((y^j)')_t = 0$.

\smallskip 
    
    These systems keep the required intersections.
    Denoting $b_j' = |S_{(\y^j)'} \setminus S_{\x^j}|$ we have
    $b_j'=b_j-1$ exactly $b_j>0$ times, and $b_j'=b_j$ exactly $n-(a_j+b_j+c_j)$ times.
    Taking the sum of (\ref{eq:sum-abc}) over all $t$,
    for the term belonging to $j$ we have
    $$
      b_j \cdot \frac{a_j! \, (b_j-1)! \, c_j!}{(a_j+(b_j-1)+c_j)!} +
      (n-a_j-b_j-c_j) \cdot \frac{a_j! \, b_j! \, c_j!}{(a_j+b_j+c_j)!}
      = n \cdot \frac{a_j! \, b_j! \, c_j!}{(a_j+b_j+c_j)!} \,,
    $$
        hence the overall sum for all $j$ is $n$ times the left-hand side of (\ref{eq:sum-abc}).
    Certainly the right-hand side is also multiplied by $n$, and the inequality follows by induction.
    
    This step is applicable unless $b_j=0$ holds for some $j$.
    Hence from now on assume $S_{\y^j}\subseteq S_{\x^j}$.

\medskip 
    
    \item As long as $S_{\x^j}\not\subseteq S_{\y^j}$ holds for all $j$, also including $j=i$:
    
    For every $t$ consider the systems
    $$
      \{ ((\x^j)', \y^j) \mid 1\leq j\leq m, \ t \notin S_{\y^j}\},
    $$
    
    where $((x^j)')_i = (x^j)_i$ for all $i \in [n]\setminus\{t\}$ and $((x^j)')_t = 0$.
    
   The argument analogous to the previous case yields the required inequality unless $a_j=0$ holds for some $j$.
    However, then we have $a_j=b_i=0$ which implies $j=i$.

\bigskip 
    
Hence for the rest of the proof assume $S_{\x^1}=S_{\y^1}$, as we can choose $i=1$, without loss of generality.
    Recall that in this situation $(x^1)_i=(y^1)_i$ for all $i \in S_{\x^1}=S_{\y^1}$.

\medskip 
    
    \item If we omit $(\x^1,\y^1)$ from the system, the left-hand side of (\ref{eq:sum-abc}) decreases by exactly 1, as currently $c_1=|S_{\x^1}|$ and $a_1=b_1=0$.
    For every $j \neq 1$ in the remaining subsystem we have $a_j,b_j>0$ because each $\x^j$ needs an entry of 2 to intersect $\y^1$, and each $\y^j$ needs an entry of 2 to intersect $\x^1$, while those two elements cannot be the same as $\x^j$ must not sum-intersect $\y^j$.

Consequently when we repeat Step (2) and (3) for the remaining system, once the procedure halts, the elements of
    $S_{\x^j} \setminus S_{\y^j}$ and of $S_{\y^j} \setminus S_{\x^j}$
    will not remain there, i.e.\ the value of the corresponding $c_j$ will be at most $\min(a,b)-1$ when $a_j=b_j=0$.

\medskip 
    
    \item The last halt occurs when the system contains a single vector-pair $(\x^j,\y^j)$ with $a_j=b_j=0$ and $c_j\geq 1$.
    This situation is reached after performing the above procedure at most $\min(a,b)-c_j+1\leq\min(a,b)$ times.
    Note that if $c_j=1$ then the intersection conditions exclude the presence of any other vector-pair.
\end{enumerate}
\end{proof}

Let us repeat that $m(k)$ denotes the maximum number of vector pairs in a strong 3-sum $(k,k)$-system and
 $$
   f(k) := \max \frac{(x+y+z)!}{x! \, y! \, z!},
 $$
 where the maximum is taken over all nonnegative integers $x,y,z$ such that $x+z\leq k$ and $y+z\leq k$. Now we are ready to prove

\begin{reptheorem}{q2s3isps-uj}

For every $k\geq 1$ we have
 $$
   f(k) \leq m(k) \leq k \cdot f(k) \,.
 $$

\end{reptheorem}

\begin{proof}[Proof of Theorem \ref{q2s3isps-uj}]
The upper bound is a consequence of Theorem \ref{t:sum-abc} as all
terms on the left-hand side of (\ref{eq:sum-abc}) are at least $(f(k))^{-1}$.
To obtain the lower bound we choose $x,y,z$ for which $f(k)$ is
attained, and choose $a=x+z$, $b=y+z$ and $c=z$ in Construction \ref{kon:abc}.
\end{proof}

\subsection{Upper bound for weak $(q+1)$-sum IVP-systems in $\{0,1,\ldots,q\}^{<\mathbb{N}}$}

Let $\{(\x^j,\y^j) \mid 1\leq j\leq m\}$ be a weak $(q+1)$-sum IVP-system in $\{0,1,\ldots,q\}^{<\mathbb{N}}$. 



\begin{obs}\label{obsi}

\

(i) For any weak $(q+1)$-sum $(a,b)$-system $\cF$ there exists another one $\cF'$  with $|\cF|=|\cF'|$ such that for any $(\x^j,\y^j)\in \cF'$ and $i$ with $x^j_i+y^j_i>0$ we have $x^j_i+y^j_i=q$.

(ii) For any strong $(q+1)$-sum $(a,b)$-system $\cF$ there exists another one $\cF'$  with $|\cF|=|\cF'|$ such that for any $(\x^j,\y^j)\in \cF'$ and $i$ with $x^j_i+y^j_i>0$ we have $x^j_i+y^j_i=q$.
\end{obs}

\begin{proof}
As $|\x^j \cap_{q+1} \y^j|=0$ implies $x^j_i+y^j_i\le q$, and increasing a coordinate helps to intersect other vectors, we can replace $\y^j$ by $\y^{{j}'}$ with $y^{{j}'}_i =q-x^j_i$.
\end{proof}

We will say that a weak/strong $(q+1)$-sum $(k,k)$-system is \textit{saturated} if it satisfies the property of Observation \ref{obsi}. For such $\cF=\{(\x^j,\y^j):1\le j \le m\}$, let us write $A^j_i$ to denote $\{t:x^j_t=i\}$ and $\alpha^j_i$ to denote $|A^j_i|$.





\begin{theorem}\label{random}
Let $p_i$ for $i=0,1,\dots,q$ be non-negative reals with $\sum_{i=0}^qp_i=1.$ If $\cF=\{(\x^j,\y^j):1\le j \le m\}$ is a saturated weak $(q+1)$-sum IVP-system, then $\sum_{j=1}^m\prod_{i=0}^qp_i^{\alpha^j_i}\le 1$ holds.
\end{theorem}


\begin{proof}
Let $(X_0,X_1,\dots,X_q)$ be a partition of $[n]$ taken at random by the rule
 $$
   \mathbb{P}(t \in X_0) = p_0 \,, \quad \mathbb{P}(t\in X_1) = p_1 \,, \quad \dots ,\quad \mathbb{P}(t \in X_q) = p_q \,,
 $$
  applied independently for each $t \in [n]=:\bigcup_{j=1}^m(S(\x^j)\cup S(\y^j))$.
For $j=1,\dots,m$ consider the events
 $$
   E_j = \bigwedge_{i=0}^q(A^j_i \subseteq X_i) \ .
 $$
We then have
 $$
   \mathbb{P}(E_j) = \prod_{i=0}^qp_i^{\alpha^j_i}\,.
 $$
Observe that $\mathbb{P}(E_j\land E_{j'})=0$ holds for all $1 \le j\neq j' \le m$.
Indeed, otherwise $A^j_i,A^{j'}_i\subseteq X_i$ holds for all $i=0,1,\dots,q$. But then for all $i=0,1,\dots,q$ we have that $z \in X_i$ and $\x^j_z=i$ implies $\y^j_z=q-i$ or $\y^j_z=0$ and similarly $\x^{j'}_z=i$ implies $\y^{j'}_z=q-i$ or $\y^{j'}_z=0$. 
So $|\x^j\cap_{q+1}\y^{j9'}|=|\x^{j'}\cap_{q+1}\y^j|=0$. This is a contradiction as the vectors are elements of a weak $(q+1)$-sum IVP-system.

Consequently the events $E_1,\dots,E_m$ mutually exclude each other, which implies that the sum of their probabilities is at most 1.
\end{proof}

Now we prove

\begin{reptheorem}{weaklimit-uj}

For any $q\ge 1$ let $m(q,k)$ and $m'(q,k)$ denote the maximum size of a strong / weak $(q+1)$-sum $(k,k)$-system. Then $\lim_{k\rightarrow \infty}\sqrt[k]{m(q,k)}=\lim_{k\rightarrow \infty}\sqrt[k]{m'(q,k)}=(\sqrt{q}+1)^2$.

\end{reptheorem}




\begin{proof}
Let us prove the upper bound first. By Observation \ref{obsi}, we can assume that $\cF$ is saturated. Then we apply Theorem \ref{random} with $p_0=p_q=\frac{\sqrt{q}-1}{q-1}$ and $p_1=p_2=\dots =p_{q-1}=p_0^2=\frac{q+1-2\sqrt{q}}{(q-1)^2}$. (Observe that $2p_0+(q-1)p_0^2=1$ as required.) As the system is saturated we have $\alpha^j_0= k-\sum_{i=1}^{q-1}\alpha^j_i=\alpha^j_q$ and thus we obtain 
$$
\prod_{i=0}^qp_i^{\alpha^j_i}= p_0^{\alpha^j_0+\alpha^j_q}p_0^{2(k-\sum_{i=1}^{q-1}\alpha^j_i)}=p_0^{2k} \ .
$$
Therefore, Theorem  \ref{random} implies $|\cF|\le (p_0^{-2})^k=((\frac{q-1}{\sqrt{q}-1})^2)^k=(\sqrt{q}+1)^{2k}$.

The lower bound is obtained using Construction \ref{kon:abc}. For fixed $q$ and growing $N$, we let $\alpha_i=p_iN$ for $i=0,1,\dots,q$ with $p_i$ as above in the proof of the lower bound, and so $k=\frac{p_0+(q-1)p_0^2}{2p_0+(q-1)p_0^2}N=(p_0+(q-1)p_0^2)N$. Then the number of pairs in the construction is $\prod_{i=0}^q\binom{(1-\sum_{j=0}^{i-1}p_j)N}{p_iN}$. Using Stirling's formula and omitting polynomial terms, this is
\[
\left[\frac{1}{p_0^{2p_0}(p_0^2)^{(q-1)p_0^2}}\right]^N=(p_0^{-2})^{(p_0+(q-1)p_0^2)N}=(p_0^{-2})^k=\left(\frac{q-1}{\sqrt{q}-1}\right)^{2k}=(\sqrt{q}+1)^{2k}.
\]
Taking $k$th root yields the claimed lower bound.
\end{proof}

\section{Concluding remarks}

There exist lots of intersection theorems all waiting to be addressed in the sum-intersection setting. We just would like to point out one. Katona's intersection theorem \cite{K64} gives the maximum size of a non-uniform $t$-intersecting family $\cF\subseteq 2^{[n]}$. 
The extremal family consists of all sets of size at least $\frac{n+t}{2}$ if $n+t$ is even, while if $n+t$ is odd, then the extremal family consists of all sets of size at least $\lceil\frac{n+t}{2}\rceil$ together with  $\binom{[n-1]}{\lfloor \frac{n+t}{2}\rfloor}$. 
One would hope to see a similar result for non-uniform $s$-sum $t$-intersecting families. That is extremal families are expected to consist of vectors of large rank. This is not going to hold as for two such vectors $\x,\y$ there might be coordinates where they $s$-sum `intersect very much' (i.e.\ $x_i+y_i$ is much larger than~$s$), but do not intersect anywhere else, so it is not a must that the support of the vectors are large.

To remedy this situation, we can define the size of the \textit{multi}-$s$-sum intersection of two vectors $\x,\y\in Q^n$ as $|\x\cap_{m,s} \y|=\sum_{i=1}^n(x_i+y_i-s+1)^+$, where for any real $z$ we define $z^+:=\max\{0,z\}$. A family $\cF\subseteq Q^n$ is \textit{$s$-multisum $t$-intersecting} if for any $\x,\y\in \cF$ we have $|\x\cap_{m,s}\y|\ge t$. Below, we  show the first step towards such intersection theorems. Katona's tool was his intersecting shadow theorem and we will need a similar result. 

We need to define the well-known shifting operation $\tau_{i,j}$ for the vector setting. For a vector $\x$ of length $n$ and integers $1\le i\neq j \le n$ we let $\tau_{i,j}(\x)$ be the vector obtained from $\x$ by exchanging its $i$th and $j$th coordinates if $x_i<x_j$ and we let $\tau_{i,j}(\x)=\x$ otherwise. For a family $\cF$ of vectors we define $\tau_{i,j}(\cF)=\{\tau_{i,j}(\x): \x\in \cF, \tau_{i,j}(\x)\notin \cF\}\cup \{\x\in \cF: \tau_{i,j}(\x)\in \cF\}$. 

The next lemma shows two basic properties of the shifting operation that are well-known for set systems.

\begin{lemma}\label{shift}
For any $\cF\subseteq Q(n,r)$ and $1\le i,j\le n$ we have $|\Delta(\tau_{i,j}(\cF))|\le |\Delta(\cF)|$. Furthermore, if $\cF$ is $s$-multisum $t$-intersecting, then so is $\tau_{i,j}(\cF)$.
\end{lemma}

\begin{proof}
Let us start with the proof of the claim concerning $t$-intersection. Suppose for $\x,\y \in\tau_{i,j}(\cF)$ we have $|\x\cap_{m,s}\y|<t$. We cannot have $\x,\y\in \cF$, as it is impossible by the $s$-multisum $t$-intersecting property of $\cF$. If  $\x,\y\in \tau_{i,j}(\cF)\setminus\cF$, then $\tau_{j,i}(\x),\tau_{j,i}(\y)\in \cF$ and $t>|\x\cap_{m,s}\y|=|\tau_{j,i}(\x)\cap_{m,s}\tau_{j,i}(\y)|$ contradicts the $s$-multisum $t$-intersecting property of $\cF$. Finally, if $\x\in \cF$ and $\y\in \tau_{i,j}(\cF)\setminus \cF$, then $\y':=\tau_{j,i}(\y)\in \cF\setminus \tau_{i,j}(\cF)$. So if $\x=\tau_{i,j}(\x)$, then $t>|\x \cap_{m,s} \y|=|\x\cap_{m,s}\y'|$ contradicts the $s$-multisum $t$-intersecting property of $\cF$. If $\x':=\tau_{i,j}(\x)\neq \x$, then as $\x\in \tau_{i,j}(\cF)$, we must have $\x'\in \cF$, and thus $t>|\x\cap_{m,s}\y|=|\x'\cap_{m,s}\y'|$ contradicts the $s$-multisum $t$-intersecting property of $\cF$. This finishes the proof that shifting preserves multisum intersecting properties.

To see $|\Delta(\tau_{i,j}(\cF))|\le |\Delta(\cF)|$ we define an injection $\iota :\Delta(\tau_{i,j}(\cF))\setminus \Delta(\cF) \rightarrow \Delta(\cF) \setminus \Delta(\tau_{i,j}(\cF))$ by letting $\iota(\x)$ be the vector obtained from $\x$ by interchanging its $i$th and $j$th coordinate. 
This is clearly an injection, all we need to verify is that every image belongs to $\Delta(\cF) \setminus \Delta(\tau_{i,j}(\cF))$. 
So let $\x\in \Delta(\tau_{i,j}(\cF))\setminus \Delta(\cF)$ be arbitrary. 
Then there exists $\y \in \tau_{i,j}(\cF)\setminus \cF$ with $\x\in \Delta(\y)$ and $\y':=\tau_{j,i}(\y)\in \cF\setminus \tau_{i,j}(\cF)$. 
Clearly, $\iota(\x)\in \Delta(\y')\subset \Delta(\cF)$. 
It remains to show $\iota(\x)\notin \Delta(\tau_{i,j}(\cF))$. 
First we claim $x_i>x_j$. 
Indeed, as $\y\in \tau_{i,j}(\cF)\setminus \cF$, we have $y_i>y_j$ showing $x_i\ge x_j$, and $x_i=x_j$ would mean $\x\in \Delta(\y')$ and $\x\in \Delta(\cF)$ contradicting $\x \in \Delta(\tau_{i,j}(\cF))\setminus \Delta(\cF)$. 
Now, $x_i>x_j$ implies $\iota(\x)_i<\iota(\x)_j$. 
Assume for a contradiction that there exists $\y^*\in \tau_{i,j}(\cF)$ with $\iota(\x)\in \Delta(\y^*)$. 
Then we must have $y^*_i\le y^*_j$. This is only possible if $\tau_{i,j}(\y^*)\in \cF$. But then $\x\in \Delta(\tau_{i,j}(\y^*))\subset \Delta(\cF)$ contradicting $\x\in \Delta(\tau_{i,j}(\cF))\setminus \Delta(\cF)$. This finishes the proof.
\end{proof}

Note that Lemma \ref{shift} is not valid for $s$-sum $t$-intersection instead of $s$-multisum $t$-intersec\-tion in the case of general $t$ as, say, the family $\{(3,2),(1,3)\}$ is 4-sum 2-intersecting, while its $(1,2)$-shift $\{(3,2),(3,1)\}$ is only 4-sum 1-intersecting.

We say that $\cF$ is \textit{left-shifted} if $\tau_{i,j}(\cF)=\cF$ for all $i<j$. Whenever $\tau_{i,j}(\cF)\neq \cF$ for some $i<j$, then $w(\cF)=\sum_{\x\in \cF}\sum_{i=1}^nix_i$ strictly decreases, so starting from any family $\cF$, after a finite number of shift operations one obtains a left-shifted family. Furthermore, by Lemma \ref{shift}, the size of the shadow does not increase and intersection properties are preserved. Therefore, when proving a lower bound on the size of shadows, one can assume that $\cF$ is left-shifted. We use the notation $2(n,r)$ for the set of vectors of rank $r$ in $\{0,1,2\}^n$.

\begin{theorem}\label{intshadow}
If $\cF\subseteq 2(n,r)$ is $3$-sum intersecting, then $|\Delta(\cF)|\ge |\cF|$.
\end{theorem}

\begin{proof}
We proceed by induction on $n$. If $n<r$, then for any $\x\in 2(n,r)$ we have $|\{i: x_i=2\}|>|\{j:x_j=0\}|$. 
Therefore for any $\cF\subseteq 2(n,r)$ in the auxiliary bipartite graph $B$ with parts $\cF$ and $\Delta(\cF)$ and edges between pairs $\y\in \Delta(\x)$, we have that the degree of any vector $\x$ in $\cF$ is at least as large as the degree of any of its neighbors $\y \in \Delta(\x)$. Consequently, $|\Delta(\cF)|\ge |\cF|$ as claimed.

If $n\ge r$, then, by Lemma \ref{shift}, we can assume that $\cF$ is left-shifted. For $a=0,1,2$ we introduce $\cF_a:=\{\x\in \cF: x_n=a\}$ and $\cF_a^-:=\{\x^-:\x\in \cF_a\}$, where $\x^-$ is the vector obtained from $\x$ by omitting its last coordinate. Observe that if $\y \in \Delta(\cF_a^-)$, then $\y^{+a} \in \Delta(\cF)$, where $\y^{+a}$ is the vector obtained from $\y$ by concatenating $a$ as a last coordinate. So, by induction, $|\Delta(\cF)|\ge \sum_{a=0}^2|\Delta(\cF_a^-)|\ge \sum_{a=0}^2|\cF_a^-|=|\cF|$ if we can prove for the second inequality that $\cF_a^-$ is 3-sum intersecting for all $a=0,1,2$. This is clear for $a=0,1$ as vectors in $\cF_a$ 3-sum intersect but as their last coordinate is 0 or 1, they must 3-sum intersect among the first $n-1$ coordinates.

Finally, consider $\cF_2^-$. Suppose for a contradiction that  $\x^-,\y^-\in \cF_2^-$ with $|\x^-~\cap_3~\y^-|$ $=0$. 
If for some $i\in [n-1]$ we have $x_i=0$ and $y_i\le 1$, then $|\tau_{i,n}(\x)\cap_3 \y|=0$ holds, contradicting the 3-sum intersecting property of $\cF_2$. We derive the same contradiction if $x_i\le 1$ and $y_i=0$. But $\x^-,\y^-\in 2(n-1,r-2)$, so there are at most $r-2$ coordinates from $i\in [n-1]$ with $x_i,y_i\ge 1$, hence there exists at least one coordinate $i$ for which we get the desired contradiction. 
\end{proof}

\begin{theorem}\label{multi1}
If $\cF\subseteq 2^n$ is 3-multisum $2$-intersecting, then $|\cF|\le |\cup_{r=n+1}^{2n}2(n,r)|$.
\end{theorem}

\begin{proof}
Let $\cF$ be a 3-multisum 2-intersecting family of maximum size. Clearly, $\cF$ is upward closed, i.e.\ $\y>\x\in\cF$ implies $\y\in\cF$. Observe that writing $\nabla(\x)=\{\y>\x:r(\y)=r(\x)+1\}$, we have that for any $\x\in \cF$ the shade $\nabla(\overline{\x})$ is disjoint from $\cF$. Let $r$ be the rank of a smallest ranked vector in $\cF$ and consider $\cF_r=\{\x\in \cF: r(\x)=r\}$. Observe that $\cF':=(\cF\setminus \cF_r) \cup \nabla(\overline{\cF_r})$ is 3-multisum 2-intersecting. Indeed, vectors from $\cF'\setminus \cF$ are all of rank $2n-r+1$ and vectors from $\cF\cap \cF'$ are all of rank at least $r+1$, so they must 3-multisum 2-intersect. As $|\nabla(\overline{\cF})|=|\Delta(\cF)|$, by Theorem \ref{intshadow}, $|\cF'|\ge |\cF|$ and we can repeat this procedure as long as $r\le n$ and thus $2n-r+1>r$. We obtain that $|\cF|\le |\cup_{r=n+1}^{2n}2(n,r)|$.
\end{proof}

\begin{theorem}\label{multi2}
If $\cF\subseteq 2^n$ is 3-multisum $3$-intersecting, then $|\cF|\le |\cup_{r=n+2}^{2n}2(n,r)|+M(n)$, where $M(n)$ denotes the maximum size of a 3-multisum 3-intersecting family in $2(n,n+1)$.
\end{theorem}

\begin{proof}
The proof is almost identical to that of Theorem \ref{multi1}. Let $\cF$ be a 3-multisum 3-intersecting family of maximum size. Observe that writing $\nabla_2(\x)=\{\y>\x:r(\y)=r(\x)+2\}$, we have that for any $x\in \cF$ the 2-shade $\nabla_2(\x)$ is disjoint from $\cF$. Let $r$ be the rank of a smallest ranked vector in $\cF$ and consider $\cF_r=\{\x\in \cF: r(\x)=r\}$. Observe that $\cF':=(\cF\setminus \cF_r) \cup \nabla_2(\overline{\cF_r})$ is 3-multisum 2-intersecting. 
Indeed, vectors from $\cF'\setminus \cF$ are all of rank $2n-r+2$ and vectors from $\cF\cap \cF'$ are all of rank at least $r+1$, so they must 3-multisum 3-intersect. Note that if $\cG$ is $s$-multisum $t$-intersecting, then $\Delta(\cG)$ is $s$-multisum $(t-2)$-intersecting. 
So applying Theorem \ref{intshadow} twice and using $|\nabla_2(\overline{\cF})|=|\Delta(\Delta(\cF))|$, we obtain $|\cF'|\ge |\cF|$ and we can repeat this procedure as long as $r\le n$ and thus $2n-r+2>r$. We obtain that there exists a maximum-sized 3-multisum 2-intersecting family $\cF$ consisting only of vectors of rank at least $n+1$.
\end{proof}



\end{document}